\documentclass[12pt]{amsart}
\usepackage{amssymb}
\usepackage{amsmath}
\usepackage{amsfonts}
\usepackage{color}

\numberwithin{equation}{section}

\theoremstyle{plain}
\newtheorem{theorem}{Theorem}[section]
\newtheorem{proposition}[theorem]{Proposition}
\newtheorem{lemma}[theorem]{Lemma}
\newtheorem{corollary}[theorem]{Corollary}

\theoremstyle{remark}

\newtheorem{case}{{\sc Case}}

\newcommand{\reft}[1]{Theorem \ref{thm:#1}}
\newcommand{\refp}[1]{Proposition \ref{prop:#1}}
\newcommand{\refps}[2]{Propositions \ref{prop:#1} and \ref{prop:#2}}
\newcommand{\refl}[1]{Lemma \ref{lem:#1}}

\newcommand{\refs}[1]{Section \ref{sec:#1}}

\newcommand{\refeq}[1]{\eqref{eq:#1}}

\newcommand{\alp}{\alpha}
\newcommand{\bet}{\beta}
\newcommand{\gam}{\gamma}

\newcommand{\del}{\delta}
\newcommand{\Del}{\Delta}
\newcommand{\eps}{\varepsilon}
\newcommand{\lam}{\lambda}
\newcommand{\sig}{\sigma}
\newcommand{\tet}{\vartheta}
\newcommand{\vphi}{\varphi}
\newcommand{\zet}{\zeta}

\newcommand{\C}{\mathbb{C}}

\newcommand{\cc}{\mathcal{C}}
\newcommand{\cd}{\mathcal{D}}
\newcommand{\ce}{\mathcal{E}}
\newcommand{\ck}{\mathcal{K}}

\newcommand{\ct}{\mathcal{T}}

\newcommand{\cx}{\mathcal{X}}

\newcommand{\hf}{\hat{f}}

\newcommand{\hk}{\hat{K}}

\newcommand{\fq}{\mathbb{F}_{q}}
\newcommand{\fqx}{\fq^{\;\times}}

\newcommand{\ovl}{\overline}
\newcommand{\sset}{\subseteq}
\newcommand{\map}[3]{#1 \colon #2 \to #3}
\newcommand{\frob}[2]{\langle #1 , #2 \rangle}
\newcommand{\set}[2]{\{ #1 \colon #2 \}}

\newcommand{\seq}[2]{ #1_{1}, \ldots, #1_{#2}}
\newcommand{\x}{\times}
\newcommand{\inv}{^{-1}}
\newcommand{\Cx}{\C^{\x}}

\newcommand{\bcap}{\bigcap}
\newcommand{\bcup}{\bigcup}
\newcommand{\bca}{\begin{cases}}
\newcommand{\eca}{\end{cases}}
\newcommand{\lpar}{\left(}
\newcommand{\rpar}{\right)}

\newcommand{\frg}{{\mathfrak{g}}}

\newcommand{\fru}{{\mathfrak{u}}}

\renewcommand{\sp}{Sp_{2n}(q)}
\newcommand{\eo}{O_{2n}(q)}
\newcommand{\oo}{O_{2n+1}(q)}
\newcommand{\frsp}{\mathfrak{sp}_{2n}(q)}
\newcommand{\freo}{\mathfrak{o}_{2n}(q)}
\newcommand{\froo}{\mathfrak{o}_{2n+1}(q)}

\newcommand{\xa}{\xi_{\alpha,r}}
\newcommand{\la}{\lambda_{\alpha,r}}
\newcommand{\oa}{O_{\alpha,r}}

\newcommand{\ld}{\lambda_{D,\phi}}

\newcommand{\xd}{\xi_{D,\phi}}
\newcommand{\ed}{e_{D,\phi}}
\newcommand{\ced}{e_{\cd,\vphi}}

\newcommand{\zd}{z_{D,\phi}}
\newcommand{\zdd}{z_{D',\phi'}}
\newcommand{\od}{O_{D,\phi}}
\newcommand{\vd}{V_{D,\phi}}
\newcommand{\kd}{K_{D,\phi}}
\newcommand{\odd}{O_{D',\phi'}}
\newcommand{\vdd}{V_{D',\phi'}}
\newcommand{\kdd}{K_{D',\phi'}}

\newcommand{\pd}{\map{\phi}{D}{\fqx}}
\newcommand{\cpd}{\map{\vphi}{\cd}{\fqx}}
\newcommand{\vpd}{\map{\vphi}{\ce(D)}{\fqx}}

\newcommand{\vod}{V_{\cd,\vphi}}

\newcommand{\xdd}{\xi_{D',\phi'}}
\newcommand{\pdd}{\map{\phi'}{D'}{\fqx}}

\newcommand{\zij}{\zet_{i,j,r}}
\newcommand{\mij}{\mu_{i,j,r}}

\newcommand{\irr}{\operatorname{Irr}}

\newcommand{\scf}{\operatorname{scf}}
\newcommand{\tr}{\operatorname{Tr}}
\newcommand{\sgn}{\operatorname{sgn}}

\newcommand{\all}{\text{ for all }}

\allowdisplaybreaks

\begin{document}


\title[]{A Supercharacter Theory for the Sylow $p$-subgroups of the finite symplectic and orthogonal groups}

\author[]{Carlos A. M. Andr\'e \& Ana Margarida Neto}

\address[C. A. M. Andr\'e]{Departamento de Matem\'atica \\ Faculdade de Ci\^encias da Uni\-ver\-si\-da\-de de Lis\-boa \\ Cam\-po Grande \\ Edi\-f\'\i \-cio C6 \\ Piso 2 \\ 1749-016 Lisboa \\ Portugal}

\address[A. M. Neto]{Instituto Superior de Economia e Gest\~ao \\ Universidade T\'ecnica de Lisboa \\ Rua do Quelhas 6 \\ 1200-781 Lisboa \\ Portugal}

\address[C. A. M. Andr\'e \& A. M. Neto]{Centro de Estruturas Lineares e Combinat\'orias \\ Complexo Interdisciplicar da Universidade de Lisboa \\ Av. Prof. Gama Pinto 2 \\ 1649-003 Lisboa \\ Portugal}

\email{caandre@fc.ul.pt}

\email{ananeto@iseg.utl.pt}


\date{\today}


\thanks{This research was made within the activities of the Centro de Estruturas Lineares e Combinat\'orias (University of Lisbon, Portugal) and was partially supported by the Funda\c c\~ao para a Ci\^encia e Tecnologia (Lisbon, Portugal) through the project POCTI-ISFL-1-1431. A large part of the research of the first author was made and concluded while he was visiting the University of Stanford (USA) and participating in the special program on ``Combinatorial Representation Theory'' at the MSRI (Berkeley, USA) whose hospitality is gratefully acknowledged, and was partially supported by the sabbatical research grant 4/2008 of the Funda\c c\~ao Luso-Americana para o Desenvolvimento (Lisbon, Portugal). The first author also expresses his sincere gratitude to Persi Diaconis for his invitation to visit the University of Stanford, and for many enlightening discussions regarding supercharacters and their applications.}

\begin{abstract}
We define the superclasses for a classical finite unipotent group $U$ of type $B_{n}(q)$, $C_{n}(q)$, or $D_{n}(q)$, and show that, together with the supercharacters defined in \cite{AN2}, they form a supercharacter theory in the sense of \cite{DI}. In particular, we prove that the supercharacters take a constant value on each superclass, and evaluate this value. As a consequence, we obtain a factorization of any superclass as a product of elementary superclasses. In addition, we also define the space of superclass functions, and prove that it is spanned by the supercharacters. As as consequence, we (re)obtain the decomposition of the regular character as an orthogonal linear combination of supercharacters. Finally, we define the supercharacter table of $U$, and prove various orthogonality relations for supercharacters (similar to the well-known orthogonality relations for irreducible characters).
\end{abstract}

\maketitle


\section{Introduction} \label{sec:intro}

This paper is a continuation of the authors' papers \cite{AN1,AN2}, and develops a supercharacter theory for the Sylow  $p$-subgroup of the symplectic or orthogonal groups defined over the finite field $\fq$ with $q$ elements; throughout the paper, $p$ will stand for an arbitrary odd prime number, and $q = p^{e}$, $e \geq 1$, will be a fixed power of $p$. The concept of a ``supercharacter theory'' for an arbitrary finite group was developed by P. Diaconis and I. M. Isaacs in the paper \cite{DI}. Roughly, a supercharacter theory replaces irreducible characters by ``supercharacters'', and conjugacy classes by ``superclasses'', in such a way that a ``supercharacter table'' can be constructed as an ``almost unitary'' matrix with similar properties as the usual character table (namely, orthogonality of rows and columns). More precisely, given any finite group $G$, a {\it supercharacter theory} for $G$ consists of  a partition $\ck$ of $G$ and a set $\cx$ of (complex) characters of $G$ satisfying the following three axioms:
\begin{enumerate}
\item[(S1)] $|\ck| = |\cx|$;
\item[(S2)] every irreducible character of $G$ is a constituent of a unique $\xi \in \cx$;
\item[(S3)] the characters in $\cx$ are constant on the member of $\ck$.
\end{enumerate}
The elements of $\ck$ will be referred to as {\it superclasses}, and the elements of $\cx$ as {\it supercharacters} of $G$. (We observe that, by \cite[Lemma~2.1]{DI}, axiom (S2) is equivalent to requiring that $\{1\} \in \ck$.)

Every finite group $G$ has two ``trivial'' supercharacter theories: the full character theory (where $\cx$ consists of all irreducible characters of $G$, and $\ck$ of all the conjugacy classes of $G$), and the one where $\cx = \{1_{G}, \rho_{G}-1_{G}\}$ and $\ck$ consists of the sets $\{1\}$ and $G - \{1\}$; as usual, we denote by $1_{G}$ the trivial character and by $\rho_{G}$ the regular character of $G$. Although for some groups these are the only possibilities, there are many groups for which nontrivial supercharacter theories exist, and in many cases it may be possible to obtain useful information using some particular supercharacter theory. An illustrating example can be found in the paper \cite{ADS} where E. Arias-Castro, P. Diaconis and R. Stanley showed that a special supercharacter theory can be applied to study a random walk on upper triangular matrices over finite fields using techniques that traditionally required the knowledge of the full character theory.

Supercharacters theories were initially developed for the upper unitrangular group $U_{n}(q)$ consisting of all unipotent upper-triangular $n \times n$ matrices over the finite field $\fq$ with $q$ elements (where $q$ is a power of some prime number $p$). In his PhD thesis \cite{A0}, the first author begun the study of the ``basic characters'' of $U_{n}(q)$ (under the assumption that $p \geq n$), and was able to show that ``clumping'' together some of the conjugacy classes and some of the irreducible characters one attains a workable ``approximation'' to the representation theory of $U_{n}(q)$. His results were published in a series of papers in the Journal of Algebra, and showed in particular that the ``basic characters'' determine uniquely the superclasses of a supercharacter theory for $U_{n}(q)$. The original theory relies on a construction due to D. Kazdhan (see \cite{K}) and is based on Kirillov's method of coadjoint orbits (see \cite{Ki} for a description of Kirillov's method for the unitriangular group). Later, in his PhD thesis \cite{Y}, N. Yan showed how the ``basic characters'' can be obtained using more elementary methods which avoid Kazhan's construction and the algebraic geometry involved in it. Yan's approach is valid for an arbitrary prime, and was generalized later by P. Diaconis and M. Isaacs in the paper \cite{DI} in order to extend the theory for an arbitrary finite algebra group defined over $\fq$.

In the papers \cite{AN1} (for sufficiently large primes) and \cite{AN2} (for arbitrary odd primes), the authors started to develop a supercharacter theory for a Sylow $p$-subgroup $U$ of one the (non-twisted) Chevalley groups $C_{n}(q)$, $B_{n}(q)$, and $D_{n}(q)$, by defining the supercharacters of $U$ and proving some of their main properties. As in the case of the unitriangular group, the supercharacters of $U$ are parametrized by certain ``minimal'' subsets of (positive) roots. In fact, it is known that the supercharacters of $U_{n}(q)$ can be obtained as certain ``reduced'' products of ``elementary characters'' which are irreducible characters corresponding to the ``matrix entries'' $(i,j)$, for $1 \leq i < j \leq n$, labelled by nonzero elements of $\fq$; in Yan's thesis, the ``elementary characters'' were called ``primary characters'', and the supercharacters were called ``transition characters''. Following Yan's method, one can show that the supercharacters of $U_{n}(q)$ are parametrized by certain combinatorial data consisting of a ``basic set'' $D$ of matrix entries such that no two elements of $D$ agree in, either the first, or the second, coordinate, and of a map $\phi$ from $D$ to the nonzero elements of $\fq$. In the papers \cite{AN1,AN2}, the authors defined the supercharacters also as certain ``reduced'' products of ``elementary characters'' (which in general are not necessarily irreducible characters) of the given Sylow $p$-subgroup $U$. These ``reduced'' products are parametrized by pairs consisting of a conveniently chosen ``basic subset of roots'' and of a map to the nonzero elements of $\fq$. (We note that the roots in the unitriangular case are in one-to-one correspondence with the matrix entries.) In fact, the group $U$ can be naturally identified with a subgroup of a unitriangular group, and the supercharacters of $U$ can be obtained as constituents of the restriction of a supercharacter of that unitriangular group.

On the other hand, as shown in Yan's thesis, the same combinatorial data parametrize the superclasses of $U_{n}(q)$ (see also \cite{A2,A3} where the superclasses were called ``basic subvarieties''). In fact, a superclass of $U_{n}(q)$ can be obtained as a ``basic'' product of ``elementary superclasses'' which are conjugacy classes corresponding to the ``matrix entries'', and labelled by nonzero elements of $\fq$. As for supercharacters, we define the superclasses of $U$ by ``restricting'' the superclasses of the larger unitriangular group, and show that they are indexed by the same combinatorial data consisting of a conveniently chosen ``basic set of roots'' where each root is labeled by a nonzero element of $\fq$.

The paper is organized as follows. In \refs{supchar}, we introduce the necessary notation, and recall the definition and main property of the supercharacters of the group $U$. In \refs{supclass}, we define the superclasses of $U$ by intersecting with the superclasses of the unitriangular group which contains $U$, and show that the set $\ck$ of all superclasses gives a partition of $U$. Then, in \refs{scf}, we define a superclass function on $U$ as a function taking a constant (complex) value on each superclass, and prove that the supercharacters are superclass functions, and form an orthogonal basis for the complex vector space $\scf(U)$ consisting of all superclass functions. As a consequence, we obtain an explicit decomposition of the regular character as a linear combination with positive integers of all the supercharacters (see \cite[Theorem~5.2]{AN2} for a different proof), and (re)prove the main theorem on supercharacters which states that every irreducible character is a constituent of a (unique) supercharacter. In \refs{value}, we determine the constant value of a supercharacter on a superclass, and conclude that the superclasses and supercharacters satisfy the axioms of a supercharacter theory for $U$ in the sense of Diaconis and Isaacs. As a consequence, we show that every superclass factorizes uniquely as a product of elementary superclasses. Finally, in \refs{table}, we define the supercharacter table of $U$ as the square matrix with entries given by the supercharacter values, and prove the main orthogonality relations for supercharacters. As a consequence, we also deduce that the space $\scf(U)$ of superclass functions is a commutative semisimple algebra with respect to the convolution product.


\section{Supercharacters} \label{sec:supchar}

Let $p \geq 3$ be a prime number, $q = p^{e}$ ($e \geq 1$) a power of $p$, and $\fq$ the finite field with $q$ elements. For a fixed positive integer $n$, let $G$ denote one of the following classical finite groups: the symplectic group $\sp$, the even orthogonal group $\eo$, or the odd orthogonal group $\oo$ (in alternative notation, these are the (non-twisted) Chevalley groups $C_{n}(q)$, $B_{n}(q)$, and $D_{n}(q)$, respectively). Throughout the paper, we set $U = G \cap U_{m}(q)$ where $$m = \bca 2n, & \text{if $G = \sp$, or $G = \eo$,} \\ 2n+1, & \text{if $G = \oo$,} \eca$$ and $U_{m}(q)$ denotes the upper unitriangular group consisting of all unipotent upper-triangular $m \times m$ matrices over $\fq$. Then, $U$ is a Sylow $p$-subgroup of $G$, and it is described as follows. Let $J = J_{n}$ be the $n \times n$ matrix with 1's along the anti-diagonal and 0's elsewhere. Then, $U$ consists of all (block) matrices of the form
\begin{equation} \label{eq:e1}
\begin{pmatrix} x & xu & xz \\ 0 & I_{r} & -u^{T}J \\ 0 & 0 & Jx^{-T}J \end{pmatrix}
\end{equation}
where $x \in U_{n}(q)$, $u$ is an $n \x r$ matrix over $\fq$, and
\begin{enumerate}
\item $r = 0$, and $Jz^{T} - zJ = 0$, if $U \leq \sp$;
\item $r = 0$, and $Jz^{T} + zJ = 0$, if $U \leq \eo$;
\item $r = 1$, and $Jz^{T} + zJ = -uu^{T}$, if $U \leq \oo$.
\end{enumerate}

As mentioned in the Introduction, both supercharacters and superclasses of $U$ will be parametrized by certain subsets of (positive) roots. Thus, we introduce some notation and recall some elementary facts concerning roots; for the details, we refer to the books \cite{C1,C2} by R. Carter (see also \cite[Chapter~8]{CR}). Let $T$ be the maximal torus of $G$ consisting of all diagonal matrices, and $\Sigma$ the root system defined by $T$. The elements of $\Sigma$ are described as follows. For each $1 \leq i \leq n$, let $\map{\eps_{i}}{T}{\fqx}$ be the map defined by $\eps_{i}(t) = t_{i}$ for all $t \in T$; here, we denote by $t_{i} \in \fqx$ the $(i,i)$th entry of the matrix $t \in T$. Then, $\Sigma = \Phi \cup (-\Phi)$ where $$\Phi = \set{\eps_{i} \pm \eps_{j}}{1 \leq i < j \leq n} \cup \Phi'$$ and $$\Phi' = \bca \set{2\eps_{i}}{1 \leq i \leq n}, & \text{if $G = \sp$,} \\ \emptyset, & \text{if $G = \eo$,} \\ \set{\eps_{i}}{1 \leq i \leq n}, & \text{if $G = \oo$.} \eca$$ The roots in $\Phi$ are said to be {\it positive}, and the roots in $-\Phi$ are said to be {\it negative}. Throughout the paper, the word ``root'' will always stand for ``positive root''.

With $\Phi$ we associate the subset of ``matrix entries'' $\ce \sset \set{(i,j)}{-n \leq i, j \leq n}$ as follows. For any $\alp \in \Phi$, we set $$\ce(\alp) = \bca \{(i,j), (-j,-i)\}, & \text{if $\alp = \eps_{i} - \eps_{j}$ for $1 \leq i < j \leq n$}, \\ \{(i,-j), (j,-i)\}, & \text{if $\alp = \eps_{i} + \eps_{j}$ for $1 \leq i < j \leq n$}, \\ \{(i,-i)\}, & \text{if $G = \sp$ and $\alp = 2\eps_{i}$ for $1 \leq i \leq n$,} \\ \{(i,0), (0,-i)\}, & \text{if $G = \oo$ and $\alp = \eps_{i}$ for $1 \leq i \leq n$,} \eca$$ and we define $$\ce = \bcup_{\alp \in \Phi} \ce(\alp).$$ More generally, for each subset $\Psi \sset \Phi$, we set $$\ce(\Psi) = \bcup_{\alp \in \Psi} \ce(\alp);$$ hence, $\ce = \ce(\Phi)$. 

On the other hand, we consider the mirror order $\prec$ on the set $\{0, \pm 1, \ldots, \pm (n+1)\}$ which is defined as $$1 \prec 2 \prec \cdots \prec n+1 \prec 0 \prec -(n+1) \prec \cdots \prec -2 \prec -1,$$ and we shall index the rows (from left to right) and columns (from top to bottom) of any $m \x m$ matrix according to this ordering. Hence, the entries of any matrix $x \in U_{m}(q)$ are indexed by all the pairs $(i,j) \in \ce$: for each $(i,j) \in \ce$, we shall write $x_{i,j}$ to denote the $(i,j)$th entry of $x$ (which occurs in the $i$th row and in the $j$th column). For our purposes, it is convenient to consider the set $$\ce^{+} = \set{(i,j) \in \ce}{1 \leq i \leq n,\ i \prec j \preceq -i},$$ and extend this notation to any subset $\Psi \sset \Phi$ by setting $$\ce^{+}(\Psi) = \ce(\Psi) \cap \ce^{+}.$$ We observe that there exists a one-to-one correspondence between $\Phi$ and $\ce^{+}$.

For any $\alp \in \Phi$, we define the subgroup $U_{\alp}$ of $U$ as follows:
\begin{enumerate}
\item if $\alp = \eps_{i}-\eps_{j}$ for $1 \leq i < j \leq n$, then $$U_{\alp} = \set{x \in U}{x_{i,k} = 0,\  i < k < j};$$
\item if $\alp = \eps_{i}-\eps_{j}$ for $1 \leq i < j \leq n$, then $$U_{\alp} = \set{x \in U}{x_{i,k} = x_{j,l} = 0,\ i < k \leq n,\ j \prec l \preceq 0};$$
\item if, either $\alp = 2\eps_{i}$ for $1 \leq i \leq n$ (in the case where $U \leq \sp$), or $\alp = \eps_{i}$ for $1 \leq i \leq n$ (in the case where $U \leq \oo$), then $$U_{\alp} = \set{x \in U}{x_{i,k} = 0,\ i < k \leq n}.$$
\end{enumerate}

Let $\map{\tet}{\fq}{ \Cx}$ be a non-trivial linear character of the additive group $\fq^{\;+}$ of $\fq$ (this character will be kept fixed throughout the paper; moreover, all characters will be taken over the complex field). For any $r \in \fqx$, the mapping $x \mapsto \tet(rx_{i,j})$ defines a linear character $\map{\la}{U_{\alpha}}{ \Cx}$ of $U_{\alpha}$, and we define the {\it elementary character} $\xa$ to be the induced character $$\xa = (\la)^{U}$$ (see \cite{A1} for the corresponding definition in the case of the unitriangular group; see also \cite[Corollary~5.11]{DI} and the discussion thereon).

We next define the notion of a ``basic subset of roots''. To start with, we recall that a subset $\cd \sset \ce$ is said to be {\it basic} if it contains at  most one entry from each row and at most one root from each column; in other words, $\cd \sset \ce$ is basic if $$|\set{j}{i \prec j \preceq -1,\ (i,j) \in \cd}| \leq 1 \quad \text{and} \quad |\set{i}{1 \preceq i \prec j,\ (i,j) \in \cd}| \leq 1$$ for all $-n \leq i, j \leq n$. Then, we say that $D \sset \Phi$ is a {\it basic subset} if $\cd = \ce(D)$ is a basic subset of $\ce$. (We will always use script letters to denote basic subsets of $\ce$, in contrast to basic subsets of $\Phi$ which will be mostly denoted by italic letters.)

Given any non-empty basic subset $D \sset \Phi$ and any map $\pd$, we define the supercharacter $\xd$ to be the product $$\xd = \prod_{\alp \in D} \xi_{\alp,\phi(\alp)}.$$ For convenience, if $D$ is the empty subset of $\Phi$, we consider the empty map $\pd$, and define $\xd$ to be the unit character $1_{U}$ of $U$. Let $$U_{D} = \bcap_{\alp \in D} U_{\alp} \quad \text{and} \quad \ld = \prod_{\alp \in D} (\lam_{\alp,\phi(\alp)})_{U_{D}}.$$ Then, $\ld$ is clearly a linear character of $U_{D}$ and, by \cite[Proposition~2.2]{AN1}, the supercharacter $\xd$ can be obtained as the induced character
\begin{equation} \label{eq:e2}
\xd = (\ld)^{U}.
\end{equation}
Throughout the paper, we will refer to the pair $(D,\phi)$ as a {\it basic pair} for $U$; hence, $D \sset \Phi$ is a basic subset, and $\pd$ is a map.

The main result of \cite{AN2} is the following theorem. (Given any finite group $G$, we denote by $\irr(G)$ the set of all irreducible characters of $G$, and by $\frob{\cdot}{\cdot}$ (or by $\frob{\cdot}{\cdot}_{G}$ if necessary) the Frobenius' scalar product on the complex vector space of all class functions defined on $G$.)

\begin{theorem}[{\cite[Theorem~1.1]{AN2}}] \label{thm:t1}
Let $\chi$ be an arbitrary irreducible character of $U$. Then, $\chi$ is a constituent of a unique supercharacter of $U$; in other words, there exists a unique basic subset $D \sset \Phi$ and a unique map $\pd$ such that $\frob{\chi}{\xd} \neq 0$.
\end{theorem}

We note that this theorem establishes axiom (S2) of the definition of a supercharacter theory.


\section{Superclasses} \label{sec:supclass}

In this section, we define the superclasses of $U$. This notion depends strongly on certain ``basic subvarieties'' defined by polynomial equations on the Lie algebra $\fru$ of $U$. Let $\frg$ denote one of the following classical Lie algebras defined over $\fq$: the symplectic Lie algebra $\frsp$, the even orthogonal Lie algebra $\freo$, or the odd orthogonal Lie algebra $\froo$. Then, $\fru = \frg \cap \fru_{m}(q)$ where $\fru_{m}(q)$ denotes the upper niltriangular Lie algebra consisting of all nilpotent upper-triangular $m \times m$ matrices over $\fq$. Thus, $\fru$ consists of all (block) matrices of the form
\begin{equation} \label{eq:e3}
\begin{pmatrix} a & u & w \\ 0 & 0_{r} & -u^{T}J \\ 0 & 0 & -Ja^{T}J \end{pmatrix}
\end{equation}
where $a \in \fru_{n}(q)$, $u$ is an $n \x r$ matrix over $\fq$, and
\begin{enumerate}
\item $r = 0$, and $Jw^{T} - wJ = 0$, if $\fru \leq \frsp$;
\item $r = 0$, and $Jw^{T} + wJ = 0$, if $\fru \leq \freo$;
\item $r = 1$, and $Jw^{T} + wJ = -uu^{T}$, if $\fru \leq \froo$.
\end{enumerate}
For any $\alp \in \Phi$, we will denote by $e_{\alp}$ the matrix in $\fru$ defined as follows (as usual, $1 \leq i < j \leq n$): $$e_{\alp} = \bca e_{i,j} - e_{-j,-i}, & \text{if $\alp = \eps_{i}-\eps_{j}$,} \\ e_{i,-j} + e_{j,-i}, & \text{if $\alp = \eps_{i}+\eps_{j}$ and $\fru \leq \frsp$,} \\ e_{i,-j} - e_{j,-i}, & \text{if $\alp = \eps_{i}+\eps_{j}$ and $\fru \leq \freo$ or $\fru = \froo$,} \\ e_{i,-i}, & \text{if $\fru \leq \frsp$ and $\alp = 2\eps_{i}$,} \\ e_{i,0} - e_{0,-i}, & \text{if $\fru \leq \froo$ and $\alp = \eps_{i}$.} \eca$$ It is clear that $\set{e_{\alp}}{\alp \in \Phi}$ is an $\fq$-basis of $\fru$.

Given be an arbitrary basic pair $(D,\phi)$ for $U$ with $D$ non-empty basic subset, we define the element $$\ed = \sum_{\alp \in D} \phi(\alpha) e_\alpha \in \fru.$$ Since $\fru \sset \fru_{m}(q)$, we may consider the orbit $$\vd = U_{m}(q) \ed U_{m}(q) \sset \fru_{m}(q)$$ for the natural action of $U_{m}(q) \x U_{m}(q)$ on $\fru_{m}(q)$ given by $(x,y) \cdot u = xuy\inv$ for all $x,y \in U_{m}(q)$ and $u \in \fru_{m}(q)$. By definition (see, for example, \cite{A3}, \cite{ADS}, \cite{DI}, or \cite{Y}), the superclass of $U_{m}(q)$ which contains the element $1+\ed \in U_{m}(q)$ is the subset $$1+\vd = \set{1+a}{a \in \vd} \sset U_{m}(q),$$ and we define the {\it superclass} $\kd$ of $U$ to be the intersection $$\kd = U \cap (1+\vd) = \set{x \in U}{x-1 \in \vd}.$$ We will also consider the intersection $$\od = \vd \cap \fru,$$ but observe that it is not necessarily true that $\kd = 1+\od$; however, there exists a bijection between $\kd$ and $\od$ (see \refl{l2} below).

The following result is a clear consequence of the definition; in the case where $D$ is the empty subset of $\Phi$, we consider the empty map $\pd$, and define $\ed = 0$, so that $\vd = \{0\}$ and $\kd = \{1\}$.

\begin{lemma} \label{lem:l1}
For every basic pair $(D,\phi)$ for $U$, the superclass $\kd$ is invariant under conjugation. In particular, $\kd$ is a union of conjugacy classes.
\end{lemma}

On the other hand, by \cite[Theorem~2.1]{A3}, for any basic pairs $(D,\phi)$ and $(D',\phi')$, we have $\vd \cap \vdd \neq \emptyset$ if and only if $(D,\phi) = (D',\phi')$. In fact, by the definition, $\ce(D)$ and $\ce(D')$ are basic subsets of $\ce = \ce(\Phi)$ satisfying $\ce(D) = \ce(D')$ if and only if $D = D'$, and the elements $\ed, e_{D',\phi'} \in \fru_{m}(q)$ are equal if and only if $D = D'$ and $\phi = \phi'$. Therefore, we obtain a disjoint union $\bcup_{D,\phi} \kd$ indexed by all basic pairs $(D,\phi)$ for $U$. We next show that $U$ equals this union, thus obtaining a partition $\ck$ of $U$ which satisfies axiom (S1) for the required supercharacter theory. The proof of the following auxiliary result can be found in \cite[Lemma~2.3]{AN2}.

\begin{lemma} \label{lem:l2}
Let  $(D,\phi)$ be a basic pair for $U$. Let $$z = \begin{pmatrix} x & xv & xw \\ 0 & I_{r} & -v^{T}J \\ 0 & 0 & J x^{-T} J \end{pmatrix} \in U \qquad \text{and} \qquad a_{z} = \begin{pmatrix} u & v & w \\ 0 & 0_{r} & -v^{T}J \\ 0 & 0 & -J u^{T}J \end{pmatrix} \in \fru$$ where $x = 1+u$. Then, $z \in \kd$ if and only if $a_{z} \in \vd$. Moreover, the mapping $z \mapsto a_{z}$ defines a bijection from $U$ to $\fru$.
\end{lemma}

Now, let $z \in U$ be arbitrary, and $a_{z} \in \fru$ be as in the lemma. Then, by \cite[Theorem~2.1]{A3} (see also \cite[Appendix~A]{DI}), there exists a unique basic subset $\cd \sset \ce$ and a unique map $\cpd$ such that $u \in \vod$ where $\vod = U_{m}(q) \ced U_{m}(q)$ and $$\ced = \sum_{(i,j) \in \cd} \vphi(i,j) e_{i,j} \in \fru_{m}(q).$$ We claim that $\cd = \ce(D)$ for a uniquely determined basic subset $D \sset \Phi$, and $\ced = \ed \in \fru$ for a uniquely determined map $\pd$; thus, we must have $\phi(\alp) = \vphi(i,j)$ for all $\alp \in D$ with $(i,j) \in \ce^{+}(\alp)$. To see this, we recall that $\vod$ is the zero set of a finite family of polynomial equations which are defined as follows.

For each entry $(i,j) \in \ce$, let $\cd(i,j)$ denote the subset $$\cd(i,j)=\set{(k,l) \in \cd}{i \prec k \preceq -1,\ 1\preceq l \prec j}.$$ Let $\cd(i,j) = \{(i_{1},j_{1}), \ldots, (i_{t},j_{t})\}$ where ${j_{1}}\prec {j_2} \prec \ldots \prec {j_{t}}$, and let $\sig \in S_{t}$ be the permutation such that ${i_{\sig(1)}} \prec {i_{\sig(2)}} \prec \ldots \prec {i_{\sig(t)}}$; as usual, we denote by $S_{t}$ the symmetric group of degree $t$. Then, for any $u \in \fru_{m}(q)$, we define $\Del_{i,j}^{\cd}(u)$ to be the determinant $$\Del_{i,j}^{\cd}(u) = \begin{vmatrix}  u_{i,j_{1}} & \cdots & u_{i,j_{t}} & u_{i,j} \\ u_{i_{\sig(1)},j_{1}} & \cdots & u_{i_{\sig(1)},j_{t}}  & u_{i_{\sig(1)},j}  \\ \vdots &  & \vdots & \vdots \\ u_{i_{\sig(t)},j_{1}} & \cdots & u_{i_{\sig(t)},j_{t}}  & u_{i_{\sig(t)},j} \end{vmatrix}.$$ We note that $\Del_{i,j}^{\cd}(u) = u_{i,j}$ whenever $\cd(i,j)=\emptyset$; in particular, if $\cd$ is empty, then $\Del_{i,j}^{\cd}(u) = u_{i,j}$ for all $u \in \fru_{m}(q)$. 

By \cite[Proposition~2.3]{A3}, we know that
\begin{equation} \label{eq:e4}
\vod = \set{u \in \fru_{m}(q)}{\Del_{i,j}^{\cd}(u) = \Del_{i,j}^{\cd}(\ced) \all (i,j) \in R(\cd)}
\end{equation}
where $R(\cd) = \ce - S(\cd)$ and $$S(\cd) = \bcup_{(i,j) \in \cd} \lpar \set{(i,s)}{j \prec s \preceq -1} \cup \set{(r,j)}{1 \preceq r \prec i} \rpar;$$ we refer to the entries in $S(\cd)$ as the {\it $\cd$-singular entries}, and to those in $R(\cd)$ as the {\it $\cd$-regular entries}. It is easy to show that, for an arbitrary $\cd$-regular entry $(i,j) \in R(\cd)$, we have
\begin{equation} \label{eq:e5}
\Del_{i,j}^{\cd}(\ced) = \bca (-1)^{t} \sgn(\sig) \vphi(i,j) \prod_{s=1}^{t} \vphi(i_{s},j_{s}), & \text{if $(i,j) \in \cd$,} \\  0, & \text{if $(i,j) \notin \cd$,} \eca
\end{equation}
where $\cd(i,j) = \{(i_{1},j_{1}), \ldots, (i_{t},j_{t})\}$, ${j_{1}}\prec {j_2} \prec \ldots \prec {j_{t}}$, and $\sig \in S_{t}$ is such that ${i_{\sig(1)}} \prec {i_{\sig(2)}} \prec \ldots \prec {i_{\sig(t)}}$.

We now prove the following auxiliary result; for simplicity of writing, given any basic subset $D \sset \Phi$ and any entry $(i,j) \in \ce$, we set $D(i,j) = \cd(i,j)$ for $\cd = \ce(D)$.

\begin{lemma} \label{lem:l3}
Let $D \sset \Phi$ be a basic subset, let $u \in \fru$, and let $(i,j) \in \ce^{+}$. Then, $$\Del_{i,j}^{\ce(D)}(u) = (-1)^{r+1} \Del_{-j,-i}^{\ce(D)}(u)$$ where $$r = \bca |D(i,j)|, & \text{if, either $\fru \leq \freo$, or $\fru \leq \froo$,} \\ |D(i,j)| , & \text{if $\fru \leq \frsp$ and $j \leq n$,} \\ |D'(i,j)| - 1, & \text{if $\fru \leq \frsp$ and $-j \leq n$,} \eca$$ and $D'(i,j) = D(i,j) \cap \set{(i,j) \in \ce^{+}}{j \leq n}$.
\end{lemma}

\begin{proof}
Let  $D(i,j) = \{(i_{1},j_{1}), \ldots, (i_{t},j_{t})\}$ where $j_{1}\prec \ldots\prec j_{t} \prec j$, and let $\sig \in S_{t}$ be such that $i \prec i_{\sig(1)}\prec \ldots i_{\sig(t)} $. By the definition of $\ce(D)$, we clearly have $$D(-j,-i) = \{(-j_{1},-i_{1}), \ldots, (-j_{t},-i_{t})\}$$ where $-j \prec -j_{t}\prec \ldots\prec -j_{1}$ and $-i_{\sig(t)}\prec \ldots\prec -i_{\sig(1)} \prec -i$. Thus, $$\Del_{-j,-i}^{\ce(D)}(u) = \begin{vmatrix} u_{-j,-i_{\sig(t)}} & \cdots & u_{-j,-i_{\sig(1)}} & u_{-j,-i} \\ u_{-j_{t},-i_{\sig(t)}} & \cdots & u_{-j_{t},-i_{\sig(1)}} & u_{-j_{t},-i} \\ \vdots &  & \vdots & \vdots \\  u_{-j_{1},-i_{\sig(t)}} & \cdots & u_{-j_{1},-i_{\sig(1)}} & u_{-j_{1},-i}  \end{vmatrix}.$$

Firstly, we assume that, either $\fru \not\leq \frsp$, or $\fru \leq \frsp$ and $j \leq n$. Let $C_{j_{1}}$, $\ldots$, $C_{j_{t}}$, $C_{j}$ denote the column vectors of $\Del^{\ce(D)}_{i,j}(u)$, and $L_{-j}$, $L_{-j_{t}}$, $\ldots$, $L_{-j_{1}}$ the row vectors of $\Del^{\ce(D)}_{-j,-i}(u)$. For any $k \in \{\seq{j}{t}, j\}$, we have $L_{-k} = - C_{k}^{\;T} J$ where $J = J_{t+1}$ is the $(t+1) \x (t+1)$ matrix with $1$'s along the anti-diagonal and $0$'s elsewhere, and thus we deduce that
{\renewcommand{\arraycolsep}{0.8mm}
\begin{align*}
\Del_{-j,-i}^{\ce(D)}(u) &= (-1)^{t+1} \begin{vmatrix} u_{i_{\sig(t)}, j} & \cdots & u_{i_{\sig(1)}, j} & u_{i, j} \\ u_{i_{\sig(t)}, j_{t}} & \cdots & u_{i_{\sig(1)}, j_{t}} & u_{i, j_{t}} \\ \vdots &  & \vdots & \vdots \\ u_{i_{\sig(t)}, j_{1}} & \cdots & u_{i_{\sig(1)}, j_{1}} & u_{i, j_{1}} \end{vmatrix} = (-1)^{t+1} \begin{vmatrix} u_{i_{\sig(t)}, j} & u_{i_{\sig(t)}, j_{t}} & \cdots & u_{i_{\sig(t)}, j_{1}} \\ \vdots & \vdots & & \vdots \\ u_{i_{\sig(1)}, j} & u_{i_{\sig(1)}, j_{t}} & \cdots & u_{i_{\sig(1)}, j_{1}} \\ u_{i, j} & u_{i, j_{t}} & \cdots & u_{i, j_{1}}  \end{vmatrix} \\ &= (-1)^{t+1} \det(J)^{2} \begin{vmatrix} u_{i,j_{1}} & \cdots & u_{i,j_{t}} & u_{i,j} \\ u_{i_{\sig(1)},j_{1}} & \cdots & u_{i_{\sig(1)},j_{t}}  & u_{i_{\sig(1)},j}  \\ \vdots &  & \vdots & \vdots \\ u_{i_{\sig(t)},j_{1}} & \cdots & u_{i_{\sig(t)},j_{t}}  & u_{i_{\sig(t)},j} \end{vmatrix} = (-1)^{t+1} \Del^{\ce(D)}_{i,j}(u),
\end{align*}
as required.}

On the other hand, suppose that $\fru \leq \frsp$ and $j = -k$ for some $1 \leq k \leq n$. By setting $i_{\sig(0)} = i$ and $j_{t+1} = j$, let $0 \leq s \leq t$ and $1 \leq s' \leq t$ be such that $i_{\sig(s)} \preceq n \prec i_{\sig(s+1)}$ and $j_{s'} \preceq n \prec j_{s'+1}$. Since $j_{1}\prec j_{2} \prec \ldots \prec j_{s'} \preceq n\prec i_{\sig(s+1)} \prec \ldots i_{\sig(t)}$, we have $u_{-j_{b},-i_{\sig(a)}} = u_{i_{\sig(a)},j_{b}} = 0$ for all $s < a \leq t$ and all $1 \leq b \leq s'$, and thus 
{\renewcommand{\arraycolsep}{0.8mm} $$\Del_{-j,-i}^{\ce(D)}(u) = \begin{vmatrix} u_{-j,-i_{\sig(t)}} & \cdots & u_{-j,-i_{\sig(s+1)}}& u_{-j,-i_{\sig(s)}} & \cdots & u_{-j,-i_{\sig(1)}} & u_{-j,-i} \\ u_{-j_{t},-i_{\sig(t)}} & \cdots & u_{-j_{t},-i_{\sig(s+1)}} & u_{-j_{t},-i_{\sig(s)}} & \cdots & u_{-j_{t},-i_{\sig(1)}} & u_{-j_{t},-i} \\ \vdots & &  \vdots & \vdots &  & \vdots   & \vdots  \\ u_{-j_{s'+1},-i_{\sig(t)}} & \cdots & u_{-j_{s'+1},-i_{\sig(s+1)}} & u_{-j_{s'+1},-i_{\sig(s)}} & \cdots & u_{-j_{s'+1},-i_{\sig(1)}} & u_{-j_{s'+1},-i} \\ 0 & \cdots & 0 & u_{-j_{s'},-i_{\sig(s)}} & \cdots & u_{-j_{s'},-i_{\sig(1)}} & u_{-j_{s'},-i} \\ \vdots & &  \vdots & \vdots &  & \vdots   & \vdots \\ 0 & \cdots & 0 & u_{-j_{1},-i_{\sig(s)}} & \cdots & u_{-j_{1},-i_{\sig(1)}} & u_{-j_{1},-i}  \end{vmatrix}.$$ Since $\fru \leq \frsp$, we deduce that
\begin{align*}
\Del_{-j,-i}^{\ce(D)}(u) &= \begin{vmatrix} -u_{i_{\sig(t)},j} & \cdots & -u_{i_{\sig(s+1)},j} & u_{{\sig(s)},j} & \cdots & u_{i_{\sig(1)},j} & u_{i,j} \\ -u_{i_{\sig(t)},j_{t}} & \cdots & -u_{i_{\sig(s+1)},j_{t}} & u_{i_{\sig(s)},j_{t}} & \cdots & u_{i_{\sig(1)},j_{t}} & u_{i,j_{t}} \\ \vdots & &  \vdots & \vdots &  & \vdots   & \vdots \\ -u_{i_{\sig(t)},j_{s'+1}} & \cdots & -u_{i_{\sig(s+1)},j_{s'+1}} & u_{i_{\sig(s)},j_{s'+1}} & \cdots & u_{i_{\sig(1)},j_{s'+1}} & u_{i,j_{s'+1}} \\ 0 & \cdots & 0 & -u_{i_{\sig(s)},j_{s'}} & \cdots & -u_{i_{\sig(1)},j_{s'}} & -u_{i,j_{s'}} \\ \vdots & &  \vdots & \vdots &  & \vdots   & \vdots  \\ 0 & \cdots & 0 & -u_{i_{\sig(s)},j_{1}} & \cdots & -u_{i_{\sig(1)},j_{1}} & -u_{i,j_{1}} \end{vmatrix} \\ &= (-1)^{t-s+s'} \begin{vmatrix} u_{i_{\sig(t)},j} & \cdots & u_{i_{\sig(s+1)},j} & u_{{\sig(s)},j} & \cdots & u_{i_{\sig(1)},j} & u_{i,j} \\ u_{i_{\sig(t)},j_{t}} & \cdots & u_{i_{\sig(s+1)},j_{t}} & u_{i_{\sig(s)},j_{t}} & \cdots & u_{i_{\sig(1)},j_{t}} & u_{i,j_{t}} \\ \vdots & &  \vdots & \vdots &  & \vdots   & \vdots \\ u_{i_{\sig(t)},j_{s'+1}} & \cdots & u_{i_{\sig(s+1)},j_{s'+1}} & u_{i_{\sig(s)},j_{s'+1}} & \cdots & u_{i_{\sig(1)},j_{s'+1}} & u_{i,j_{s'+1}} \\ 0 & \cdots & 0 & u_{i_{\sig(s)},j_{s'}} & \cdots & u_{i_{\sig(1)},j_{s'}} & u_{i,j_{s'}} \\ \vdots & &  \vdots & \vdots &  & \vdots   & \vdots  \\ 0 & \cdots & 0 & u_{i_{\sig(s)},j_{1}} & \cdots & u_{i_{\sig(1)},j_{1}} & u_{i,j_{1}} \end{vmatrix}.
\end{align*}
Arguing as above (transposing and conjugating by the matrix $J = J_{t+1}$), we conclude that $$\Del_{-j,-i}^{\ce(D)}(u) = (-1)^{t-s+s'} \Del_{i,j}^{\ce(D)}(u),$$ and the result follows because $t-s+s' = |D'(i,j)|$.}
\end{proof}

We are now able to prove the following result.

\begin{proposition} \label{prop:p1}
Let $u \in \fru$ be arbitrary, and $(\cd,\vphi)$ be the (unique) basic pair for $U_{m}(q)$ such that $u \in \vod$. Then, $\ced \in \fru$; in particular, there exists a unique basic pair $(D,\phi)$ for $U$ such that $\cd = \ce(D)$ and $\ced = \ed$.
\end{proposition}

\begin{proof}
It is enough to show that, for all $(i,j) \in \cd \cap \ce^{+}$, we have $$\vphi(i,j) = (-1)^{\eps_{j}} \vphi(-j,-i)$$ where $$\eps_{j} = \bca 1, & \text{if, either $\fru \not\leq \frsp$, or $\fru \leq \frsp$ and $j \leq n$,} \\ 0, & \text{if $\fru \leq \frsp$ and $-j \leq n$.} \eca$$ To prove this, we proceed by induction on $|\cd|$. We consider the total order $\preceq$ on the set of entries $\ce$ defined as follows: for all $(i,j),(k,l) \in \ce$, we set
\begin{equation} \label{eq:ord}
(i,j)\prec (k,l) \iff \text{either $j\prec l$, or $j=l$ and $k\prec i$.}
\end{equation}
Let $(i,j) \in \cd \cap \ce^{+}$ be the smallest entry satisfying $\vphi(i,j) \neq (-1)^{\eps_{j}} \vphi(-j,-i)$ (if it exists), let $$\cd(i,j)=\{(i_{1},j_{1}),\ldots,(i_{t},j_{t})\}$$ with $j_{1}\prec\ldots\prec j_{t}$, and let $\sig \in S_{t}$ such that $i_{\sig(1)} \prec \ldots \prec i_{\sig(t)}$. Since $(r,s) \prec (i,j)$, we have $\vphi( r,s) = (-1)^{\eps_{s}} \vphi(-s,-r)$ for all $(r,s) \in \cd(i,j) \cap \ce^{+}$. Using \refeq{e3} and \refeq{e4}, it is easy to conclude that $$\cd(-j,-i) = \{(-j_{1},-i_{1}),\ldots,(-j_{t},-i_{t})\}.$$ Let $\cd_{0} = \cd(i,j) \cup \cd(-j,-i)$, and $\map{\vphi_{0}}{\cd_{0}}{\fqx}$ be the restriction of $\vphi$ to $\cd_{0}$. Then, by induction, the element $$e_{\cd_{0},\vphi_{0}} = \sum_{(r,s) \in \cd_{0}} \vphi(r,s) e_{r,s}$$ lies in $\fru$.

On the other hand, let $c_{-j,-i} \in \fq$ denote the $(i,j)$th coefficient of $\ced$; hence, $$c_{-j,-i} = \bca \vphi(-j,-i),  & \text{if $(-j,-i) \in \cd$,} \\ 0, & \text{if $(-j,-i) \notin \cd$.} \eca$$ By \refeq{e4}, we have $$\Del_{-j,-i}^{\cd}(\ced) = (-1)^{t} \sgn(\sig) c_{-j,-i} \prod_{s=1}^{t} \vphi(-j_{s},-i_{s}).$$ It is easy to show that $$\prod_{s=1}^{t} \vphi(-j_{s},-i_{s})=(-1)^{r'}\prod_{s=1}^{t} \vphi(i_{s},j_{s})$$ where $$r' = \bca |\cd(i,j)|, & \text{if, either $\fru \leq \freo$, or $\fru \leq \froo$,} \\ |\cd'(i,j)|, & \text{if $\fru \leq \frsp$,} \eca$$ and $\cd'(i,j) = \cd(i,j) \cap \set{(i,j) \in \ce^{+}}{j \leq n}$. Thus, we deduce that $$\Del_{-j,-i}^{\cd}(\ced) = (-1)^{r'} c_{-j,-i}\, \vphi(i,j)\inv \Del_{i,j}^{\cd}(\ced).$$

Since $u \in \vod$, we have $$\Del_{i,j}^{\cd}(u) = \Del_{i,j}^{\cd}(\ced) = \Del_{i,j}^{\cd_{0}}(\ced)$$ and $$\Del_{-j,-i}^{\cd}(u) = \Del_{-j,-i}^{\cd}(\ced) = \Del_{-j,-i}^{\cd_{0}}(\ced);$$ we note that $(-j,-i) \in R(\cd)$ (by induction and by the choice of $(i,j)$). By the previous lemma, we conclude that  $$\Del_{i,j}^{\cd}(u) = (-1)^{r+1} \Del_{-j,-i}^{\cd}(\ced)$$ where $$r = \bca r'-1 , & \text{if $\fru \leq \frsp$ and $-j \leq n$,} \\ r', & \text{otherwise.} \eca$$ It follows that $$(-1)^{r'} c_{-j,-i} \vphi(i,j)\inv \Del_{i,j}^{\cd}(\ced) = (-1)^{r'+1} \Del_{i,j}^{\cd}(\ced).$$ Since $\Del_{i,j}^{\cd}(u) \neq 0$ (because $(i,j)\in \cd$), we obtain $$c_{-j,-i} = \bca \vphi(i,j), & \text{if $\fru \leq \frsp$ and $-j \leq n$,} \\ -\vphi(i,j), & \text{otherwise.} \eca$$ It follows that $(-j,-i) \in \cd$ and that $c_{-j,-i} = \vphi(-j,-i)$. Moreover, we note that, in the orthogonal case, if $j = i$, we obtain $\vphi(i,-i) = -\vphi(i,-i)$, hence $(i,-i) \notin \cd$.

The above contradicts the minimal choice of $(i,j)$, and thus we conclude that $$\vphi(i,j) = (-1)^{\eps_{j}} \phi_{\eps}(-j,-i)$$ for all $(i,j) \in \cd \cap \ce^{+}$. The proof is complete.
\end{proof}

As observed above, this concludes the proof of the following result; we recall that $\od = \fru \cap \vd$ (by the definition).

\begin{theorem} \label{thm:t2}
Let $u \in \fru$ be arbitrary. Then, there exists a unique basic subset $D \sset \Phi$ and a unique map $\pd$ such that $u \in \od$. Thus, $\fru$ is the disjoint union $$\fru = \bcup_{D,\phi} \od$$ where the union runs over all basic pairs $(D,\phi)$ for $U$. Moreover, we have $$\od = \set{a \in \fru}{\Del_{i,j}^{D}(a) = \Del_{i,j}^{D}(\ed) \all (i,j) \in R(D)}$$ for every basic pair $(D,\phi)$ for $U$. 
\end{theorem}

As a consequence of \reft{t2} and \refl{l2}, we obtain the following main theorem.

\begin{theorem} \label{thm:t3}
Let $z \in U$ be arbitrary. Then, there exists a unique basic subset $D \sset \Phi$ and a unique map $\pd$ such that $z \in \kd$. Thus, $U$ is the disjoint union of all its superclasses; that is, $$U = \bcup_{D,\phi} \kd$$ where the union runs over all basic pairs $(D,\phi)$ for $U$.
\end{theorem}


\section{Superclass functions} \label{sec:scf}

In this section, we prove that every supercharacter is a ``superclass function'' of $U$; by definition, a function $\map{\eta}{U}{\Cx}$ is said to be a {\it superclass function} if it takes a constant value on each superclass of $U$. In fact, we shall prove the following result (cf. \cite[Theorem~3.1]{AN2}).

\begin{theorem} \label{thm:t4}
Every supercharacter of $U$ is a superclass function. Moreover, every superclass function on $U$ is a linear combination of supercharacters; hence, the supercharacters of $U$ form a basis for the complex vector space $\scf(U)$ consisting of all superclass functions on $U$.
\end{theorem}

Since every supercharacter is a product of elementary characters (by definition), it is enough to show that every elementary character of $U$ takes a constant value on each superclass of $U$. In fact, the theorem above will follow from the following result. (Henceforth, for each $z \in U$, we denote by $a_{z}$ the element of $\fru$ given by \refl{l2}.)

\begin{proposition} \label{prop:p2}
Let $\alp \in \Phi$ and $r \in \fqx$ be arbitrary. Let $(D,\phi)$ be a basic pair for $U$, and denote by $\zd$ the unique element of $z \in U$ with $a_{z} = \ed$. Then, $$\xa(z) = \xa(\zd)$$ for all $z \in \kd$.
\end{proposition}

\begin{proof}[Proof (first part)]
Let $(i,j) \in \ce^{+}(\alp)$; hence, $1 \leq i \leq n$ and $i \prec j \preceq -i$. In this first part of the proof, we shall assume that $j \neq -i$ (in the case where $U \leq \sp$). Let $\zij$ be the elementary character of $U_{m}(q)$ associated with $(i,j)$ and $r$. We recall its definition (see \cite[Lemma~3]{A1}). We consider the subgroup $U_{i,j} = \set{x \in U_{m}(q)}{x_{i,k} = 0 \all i \prec k \prec j}$ of $U_{m}(q)$, and the linear character $\map{\mij}{U_{i,j}}{\Cx}$ defined by $\mij(x) = \tet(rx_{i,j})$ for all $x \in U_{i,j}$. Then, $\zij$ is defined to be the induced character $\zij = (\mij)^{U_{m}(q)}$.

By \cite[Proposition~3.2]{AN2}, we have $\xa = (\zij)_{U}$, and so $\xa(z) = \zij(z)$ for all $z \in U$. By \cite[Proposition~5.1]{A3} (see also \cite[Theorem~5.8]{DI}, or \cite[Theorem~2.2]{ADS}), we have $\zij(x) = \zij(1+\ed)$ for all $x \in 1+\vod$. In particular, we deduce that $$\xa(z) = \zij(z) = \zij(1+\ed) = \xa(\zd)$$ for all $z \in \kd$; we recall that $z \in \kd$ if and only if $a_{z} \in \vd$ (by \refl{l2}).
\end{proof}

In order to complete the proof of \refp{p2}, it remains to consider the case where $U \leq \sp$ and $\alp = 2\eps_{i}$ for some $1 \leq i \leq n$. In what follows, we will always assume that this is the case; moreover, the basic pair $(D,\phi)$ will be kept fixed. We prove some elementary auxiliary lemmas; we mention that similar results are valid in the general case (a proof can be found in the second's author PhD thesis \cite{N}). The proof of the first lemma is straightforward.

\begin{lemma} \label{lem:l4}
Let $u \in \od$ be arbitrary, and $(k,l) \in \ce(D)$ be the smallest entry of $\ce(D)$ (with respect to the total order $\preceq$ on $\ce$ as defined in \refeq{ord}); hence, $1 \leq k \leq n$ and $k \prec l \preceq -k$. Then, there exists $x \in U$ such that $v = xux\inv \in \od$ satisfies $v_{k',l} = v_{k,l'} = 0$ for all $1 \preceq k', l'  \preceq -1$ with $k' \neq k$ and $l' \neq l$.
\end{lemma}

As a consequence, we obtain the following result.

\begin{corollary} \label{cor:c1}
Let $\bet \in \Phi$ and $s \in \fqx$ be arbitrary. Then, $O_{\bet,s} = \set{x(se_{\bet})x\inv}{x \in U}$ is the adjoint $U$-orbit which contains $se_{\bet} \in \fru$. Moreover, $K_{\bet,s} = \set{xz_{\bet,s}x\inv}{x \in U}$ is the conjugacy class which contains the element $z_{\bet,s} = 1 + se_{\bet} \in U$. (In particular, we have $\xa(z) = \xa(z_{\bet,s})$ for all $z \in K_{\bet,s}$.)
\end{corollary}

\begin{proof}
The first assertion is an immediate consequence of the previous lemma. For the second, we note that $xz_{\bet,s}x\inv = 1+x(se_{\bet})x\inv \in U \cap (1+V_{\bet,s}) = K_{\bet,s}$. On the other hand, if $z \in K_{\bet,s}$, then $a_{z} \in O_{\bet,s}$ (by \refl{l2}), and thus the mapping $z \mapsto a_{z}$ defines a bijection from $K_{\bet,s}$ to $O_{\bet,s}$. Therefore, $|K_{\bet,s}| = |O_{\bet,s}| = |\set{1+x(se_{\bet})x\inv}{x \in U}|$, and the result follows.
\end{proof}

We observe that, in the notation of the corollary, we have $z_{\bet,s} = \zd$ for $D = \{\bet\}$ and $\pd$ defined by $\phi(\bet) = s$; hence, \refp{p2} is true whenever the basic subset $D$ has a unique element. Therefore, we will assume that $|D| > 1$, and $\xa(z) = \xa(z_{D',\phi'})$ for all $z \in K_{D',\phi'}$ and every basic pair $(D',\phi')$ with $|D'| < |D|$. \refp{p2} will then follow by induction. However, we need a concrete formula for the values of an elementary character on any superclass.

 As usual, we denote by $\fru^{\ast}$ the dual vector space of $\fru$, and let $\set{e^{\ast}_{\alp}}{\alp \in \Phi}$ be the $\fq$-basis of $\fru^{\ast}$ dual to the basis $\set{e_{\alp}}{\alp \in \Phi}$ of $\fru$; hence, $e^{\ast}_{\alp}(e_{\bet}) = \del_{\alp,\bet}$ for all $\alp, \bet \in \Phi$. For each $f \in \fru^{\ast}$, we define $$u(f) = \sum_{\bet \in \Phi} u_{\bet} e^{\ast}_{\bet} \in \fru_{2n}(q)$$ where $$u_{\bet} = \begin{cases} \frac{1}{2}\, f(e_{\bet}), & \text{if $\bet = \eps_{k} \pm \eps_{l}$ for $1 \leq k < l \leq n$,} \\ f(e_{\bet}), & \text{if $\bet = 2\eps_{k}$ for $1 \leq k \leq n$.} \end{cases}$$ It is easy to see that $f(v) = \tr(u(f)^{T} v)$ for all $v \in \fru$, and that the mapping $f \mapsto u(f)$ defines a vector space isomorphism from $\fru^{\ast}$ to $\fru$. Finally, we define the linear function $\hf \in \fru_{2n}(q)^{\ast}$ by $$\hf(v) = \tr(u(f)^{T} v)$$ for all $v \in \fru_{2n}(q)$, and set $$\oa^{\ast} = \set{f \in \fru^{\ast}}{\hf \in U_{2n}(q)(r e^{\ast}_{i,-i}) U_{2n}(q)}$$ where $(xgy)(a) = g(x\inv a y\inv)$ for all $x,y \in U_{2n}(q)$, $g \in \fru_{2n}(q)^{\ast}$ and $a \in \fru_{2n}(q)$. By \cite[Proposition~5.2]{AN2}, we know that $$\xa(z) = \frac{\xa(1)}{|\oa^{\ast}|} \sum_{f \in \oa^{\ast}} \tet_{f}(a_{z})$$ for all $z \in U$. Given any $\fq$-vector space $V$ and any linear map $f \in V^{\ast}$, we denote by $\tet_{f}$ the composite map $\map{\tet \circ f}{V}{\Cx}$; it is straightforward to check that $\tet_{f}$ is a linear character of the additive group $V^{+}$ and that $\irr(V^{+}) = \set{\tet_{f}}{f \in V^{\ast}}$.
 
For our purposes, it is convenient to describe the subset $\oa^{\ast} \sset \fru^{\ast}$ (and the elementary character $\xa$) as follows. By \cite[Corollary~5.3]{AN1}, we have $f \in \oa^{\ast}$ if and only if the following holds:
 \begin{enumerate}
\item $\hf(e_{a,-b}) = 0$ for all $(a,-b) \in \ce$ with $1 \leq a < i$ or $1 \leq b < i$;
 \item $\hf(e_{i,-i}) = r$;
 \item $\begin{vmatrix} \hf(e_{i,b}) & \hf(e_{i,-i}) \\ \hf(e_{a,b}) & \hf(e_{a,-i}) \end{vmatrix} = 0$ for all $(a,b) \in \ce$ with $1 \prec a \prec b \prec -i$. 
\end{enumerate}
On the other hand, for $u = u(f)$, we have $$\hf(e_{a,b}) = \tr(u^{T} e_{a,b}) = \tr(e_{a,b}^{\;T} u) = e_{a,b}^{\ast}(u) = u_{a,b}$$ for all $(a,b) \in \ce$. Therefore, for any $(a,b) \in \ce^{+}$, we deduce that $$\hf(e_{-b,-a}) = \bca \hf(e_{a,b}), & \text{if $(a,b) \in \ce^{+}(\eps_{a}-\eps_{b})$,} \\ -\hf(e_{a,b}), & \text{if $(a,b) \in \ce^{+}(\eps_{a}+\eps_{-b})$.} \eca$$ In particular, for $i \prec a \prec -i$, we get $$\hf(e_{a,-i}) = \bca \hf(e_{i,-a}), & \text{if $-n \preceq a \prec -i$,} \\ -\hf(e_{i,-a}), & \text{if $i \prec a \preceq n$,} \eca$$ and so the elements of $\oa^{\ast}$ can be parametrized by the set $\cc$ consisting of all functions $\map{c}{\set{a}{i \prec a \prec -i}}{\fq}$. In fact, for each $c \in \cc$, there exists a unique linear function $f_{c} \in \fru^{\ast}$ such that $\hf_{c} \in \fru_{2n}(q)^{\ast}$ satisfies $\hf_{c}(e_{i,a}) = c_{a}$ for all $i \prec a \prec -i$; here, we write $c_{a} = c(a)$ for all $i \prec a \prec -i$. This concludes the proof of the following result.

\begin{lemma} \label{lem:l5}
Let $\alp = 2\eps_{i} \in \Phi$ for some $1 \leq i \leq n$, and $r \in \fqx$. Then, in the notation as above, we have $$\oa^{\ast} = \set{f_{c}}{c \in \cc}.$$ Moreover, the mapping $c \mapsto f_{c}$ defines a bijection from $\cc$ to $\oa$, and $$\xa(z) = \frac{1}{q^{n-i}} \sum_{c \in \cc} \tet_{f_{c}}(a_{z})$$ for all $z \in U$.
\end{lemma}

Henceforth, we set $\tet_{c} = \tet_{f_{c}}$ for all $c \in \cc$, and consider the function $\map{\tet}{\fru}{\C}$ defined by $$\tet(u) = \sum_{c \in \cc} \tet_{c}(u)$$ for all $u \in \fru$. Since $\oa^{\ast} \sset \fru^{\ast}$ is $U$-invariant (for the natural action given by conjugation), we clearly have $\tet(xux\inv) = \tet(u)$ for all $u \in \fru$. Let $u \in \od$ be arbitrary, and $(k,l) \in \ce$ be the smallest entry of $\ce(D)$ (with respect to the order $\preceq$ defined in \refeq{ord}); hence, we must have $1 \leq k \leq n$ and $k \prec l \preceq -k$. Let $\bet \in \Phi$ be such that $(k,l) \in \ce(\bet)$ (hence, $\bet \in D$), and $s = \phi(\bet)$. By \refl{l3}, there exists $x \in U$ such that $v = xux\inv$ satisfies $v_{a,l} = v_{k,b} = 0$ for all $1 \preceq a, b  \preceq -1$ with $a \neq k$ and $b \neq l$, and thus $$\tet(u) = \tet(v) = \sum_{c \in \cc} \tet_{c}(se_{\bet} + w) = \sum_{c \in \cc} \tet_{c}(s e_{\bet}) \tet_{c}(w)$$ where $w = v - se_{\bet} \in \fru$. Next, we consider the relative positions of the entries $(k,l)$ and $(i,-i)$. There are four distinct cases.

\begin{case}[$k = i$ and $l \prec -i$]
In this case, $f_{c}(e_{\bet}) = \hf_{c}(e_{i,l} \pm e_{-l,-i}) = 2\hf(e_{i,l}) = 2c_{l}$ for all $c \in \cc$, and so $$\tet(u) = \sum_{c \in \cc} \tet(2c_{l}s) \tet_{c}(w).$$ On the other hand, we have $w_{\gam} = w_{a,b}$ whenever $\gam \in \Phi$ and $(a,b) \in \ce^{+}(\gam)$, and thus
\begin{equation} \label{eq:e6}
f_{c}(w) = \sum_{\gam \in \Phi} w_{\gam} f_{c}(e_{\gam}) = rw_{i,-i} + \sum_{i \prec b \prec -i} c_{b} w_{i,b} + \sum_{i \prec b \prec -i} \sum_{i \prec a \preceq -b} r\inv c_{-a} c_{b} w_{a,b}\,.
\end{equation}
In fact, for all $c \in \cc$ and all $\bet \in \Phi$, we have $$f_{c}(e_{\bet}) = \hf_{c}(e_{\bet}) = \bca r, & \text{if $\bet = 2\eps_{i}$,} \\ 2c_{b}, & \text{if $(i,b) \in \ce^{+}(\bet)$ and $b \neq -i$,} \\ 2 r\inv c_{-a} c_{b}, & \text{if $(a,b) \in \ce^{+}(\bet)$ for $i < a$,} \\ 0, & \text{otherwise.} \eca$$ Since $w_{a,l} = v_{a,l} = 0$ for all $i \prec a \preceq -i$, the coordinate $c_{l}$ of any $c \in \cc$ does not occur in the expression of $f_{c}(w)$ given by \refeq{e6}, and so $$\tet(u) = \lpar \sum_{t \in \fq} \tet(2ts) \rpar \lpar \sum_{c \in \cc_{1}} \tet_{c}(w) \rpar = \lpar \sum_{t \in \fq} \tet_{t}(2s) \rpar \lpar \sum_{c \in \cc_{1}} \tet_{c}(w) \rpar$$ where $\cc_{1} = \set{c \in \cc}{c_{l} = 0}$, and $\tet_{t}$, for $t \in \fq$, denotes the character of $\fq^{\;+}$ defined by $\tet_{t}(a) = \tet(ta)$ for all $a \in \fq$. Since $\sum_{t \in \fq} \tet_{t}$ is the regular character of $\fq^{\;+}$, we conclude that $\tet(u) = 0$. Furthermore, we note that $$\sum_{c \in \cc_{1}} \tet_{c}(w) = q\inv \sum_{c \in \cc} \tet_{c}(w)$$ (by the same reason as above).
\end{case}

\begin{case}[$k = i$ and $l = -i$]
We have $f_{c}(e_{\bet}) = f_{c}(e_{i,-i}) = r$ for all $c \in \cc$, and so $$\tet(u) = \tet(rs) \sum_{c \in \cc} \tet_{c}(w).$$
\end{case}

\begin{case}[$k < i$]
We have $f_{c}(e_{\bet}) = 0$ for all $c \in \cc$, and so $$\tet(u) = \sum_{c \in \cc} \tet_{c}(w).$$
\end{case}

\begin{case}[$k > i$]
In this case, we have $$f_{c}(e_{\bet}) = \bca r\inv c_{-k} c_{l}, & \text{if $l \neq -k$,} \\ r\inv c_{-k}^{\;2}, & \text{if $l = -k$,} \eca$$ for all $c \in \cc$.

On the one hand, suppose that $l \neq -k$. Then, the entries $c_{-k}$ and $c_{l}$ of any $c \in \cc$ do not occur in $f_{c}(w)$ (see the argument in case 1), hence we get $$\tet(u) = \lpar \sum_{t,t' \in \fq} \tet(r\inv s tt') \rpar \lpar q^{-2} \sum_{c \in \cc} \tet_{c}(w) \rpar.$$ Since $\sum_{t \in \fq} \tet_{t}$ is the regular character of $\fq^{\;+}$, we conclude that $$\sum_{t,t' \in \fq} \tet(r\inv s tt') = \sum_{t,t' \in \fq} \tet_{t}(r\inv s t') = q,$$ and thus $$\tet(u) = q\inv \sum_{c \in \cc} \tet_{c}(w).$$

On the other hand,  suppose that $l = -k$. Then, the entry $c_{-k}$ of any $c \in \cc$ do not occur in $f_{c}(w)$, hence we get $$\tet(u) = \lpar \sum_{t \in \fq} \tet(r\inv s t^{2}) \rpar \lpar q\inv \sum_{c \in \cc} \tet_{c}(w) \rpar.$$ Now, we recall that the {\it quadratic character} of $\fq$ is, by definition, the linear character $\eta$ of the multiplicative group $\fqx$ defined by $$\eta(c) = \bca 1, & \text{if $c \in (\fqx)^{2}$}, \\ -1, & \text{otherwise,} \eca$$ for all $c \in \fqx$. Moreover, given any linear character $\nu$ of $\fqx$ and any linear character $\tet$ of $\fq^{\;+}$, the {\it Gauss sum} of $\nu$ and $\tet$ is defined by $$G(\nu,\tet) = \sum_{c \in \fqx} \nu(c) \tet(c).$$ The following result is Theorem~5.3.3 of the book \cite{LN}.

\begin{theorem} \label{thm:t5}
Let $\tet$ be a non-trivial linear character of $\fq^{\;+}$, and $$h(T) = a_{2}T^{2} + a_{1}T + a_{0} \in \fq[T]$$ be a polynomial over $\fq$ with $a_{2} \neq 0$. Suppose that $q$ is odd. Then, $$\sum_{c \in \fq} \tet(h(c)) = \tet(a_{0} - a_{1}^{\,2} (4a_{2})\inv) \eta(a_{2}) G(\eta,\tet)$$ where $\eta$ is the quadratic character of $\fq$.
\end{theorem}

Applying this result to our situation (with $h(T) = r\inv s T^{2}$), we obtain $$\tet(u) = q\inv \eta(r\inv s) G(\eta,\tet) \sum_{c \in \cc} \tet_{c}(w).$$
\end{case}

It follows that, in any case, we have $$\tet(u) = c_{\bet,s} \sum_{c \in \cc} \tet_{c}(w)$$ for some constant $c_{\bet,s} \in \C$ depending only on the root $\bet \in D$ and on the value $s = \phi(\bet)$; in fact,
\begin{equation} \label{eq:e7}
c_{\bet,s} = \bca 0, & \text{if $k = i$ and $l \prec -i$,} \\ \tet(rs), & \text{if $k = i$ and $l = -i$,} \\ 1, & \text{if $k < i$,} \\ q\inv, & \text{if $k > i$ and $l \neq -k$,} \\ q\inv \eta(r\inv s) G(\eta,\tet), & \text{if $k > i$ and $l = -k$.} \eca
\end{equation}

We are now able to conclude the proof of \refp{p2}.

\begin{proof}[Proof of \refp{p2} (second part)]
Let the notation be as above, and let $z' \in U$ be such that $a_{z'} = w$. Let $D' = D - \{\bet\}$, and $\map{\phi'}{D'}{\fqx}$ be the restriction on $\phi$ to $D'$. Then, it is easy to check that $w = v - se_{\bet} \in \odd$, hence $z' \in \kdd$ (by \refl{l2}). By induction, we have $\xa(z') = \xa(\zdd)$. Since $$\xa(z') = \frac{1}{q^{n-i}} \sum_{c \in \cc} \tet_{c}(w)$$ (by \refl{l5}), we conclude that $$\xa(z) = \frac{c_{\bet,s}}{q^{n-i}} \sum_{c \in \cc} \tet_{c}(w) = c_{\bet,s} \xa(z') = c_{\bet,s} \xa(\zdd).$$ Therefore, the value $\xa(z)$ does not depend on $z \in \kd$, hence $\xa(z) = \xa(\zd)$ for all $z \in \kd$.
\end{proof}

We next proceed with the proof of \reft{t4}; a slightly different proof will be given later without reference to the results of \cite{AN2}.

\begin{proof}[Proof of \reft{t4}]
Since every supercharacter is a product of elementary characters, \refp{p2} implies that every supercharacter is a superclass function. By \cite[Theorem~4.2]{AN2}, the supercharacters are orthogonal, hence they are linearly independent functions of $\scf(U)$. Since the dimension of the vector space $\scf(U)$ equals the number of basic pairs $(D,\phi)$ for $U$, we conclude that the supercharacters form a basis of $\scf(U)$, and this completes the proof.
\end{proof}

We now observe that, since the regular character of $U$ is clearly a superclass function, \reft{t4} implies that it is a linear combination of supercharacters. In particular, we obtain the following result (and also an alternative proof of \cite[Theorem~3.2]{AN2}).

\begin{theorem} \label{thm:t8}
Every irreducible character is a constituent of a (unique) supercharacter. 
\end{theorem}

\begin{proof}
It is enough to observe that every irreducible character of $U$ is a constituent of the regular character. (The unicity follows by the orthogonality of supercharacters; see \cite[Theorem~4.2]{AN2}.)
\end{proof}

Finally, an easy calculation proves the following result (and gives an alternative proof of \cite[Theorem~5.2]{AN2}).

\begin{theorem} \label{thm:t9}
Let $\rho_{U}$ be the regular character of $U$. Then, $$\rho_{U} = \sum_{D,\phi} \frac{\xd(1)}{\frob{\xd}{\xd}}\; \xd$$ where the sum is over all basic pairs $(D,\phi)$.
\end{theorem}

\begin{proof}
Let $\rho_{U} = \sum_{D,\phi} m_{D,\phi} \xd$ where $m_{D,\phi} \in \C$ for all basic pairs $(D,\phi)$. Since supercharacters are orthogonal, we obtain $\frob{\rho_{U}}{\xd} = m_{D,\phi} \frob{\xd}{\xd}$. On the other hand, let $\irr_{D,\phi}(U)$ denote the subset of $\irr(U)$ consisting of all irreducible constituents of the supercharacter $\xd$; hence, we have a disjoint union $\irr(U) = \bcup_{D,\phi} \irr_{D,\phi}(U)$, and $$\xd = \sum_{\chi \in \irr_{D,\phi}(U)} \frob{\chi}{\xd} \chi.$$ Since $\rho_{U} = \sum_{\chi \in \irr(U)} \chi(1) \chi$, we deduce that $$m_{D,\phi} \frob{\xd}{\xd} = \frob{\rho_{U}}{\xd} = \sum_{\chi \in \irr_{D,\phi}(U)} \chi(1) \frob{\chi}{\xd} = \xd(1),$$ and the result follows.
\end{proof}


\section{Supercharacter values} \label{sec:value}

In this section, we obtain explicit formulae that allows to determine the constant value $\xd(\zdd)$ of the supercharacter $\xd$ on the superclass $\kdd$. Since $\xd$ is a product of elementary characters, it is enough to determine the value of an arbitrary elementary character $\xa$, for $\alp \in \Phi$ and $r \in \fqx$, on any superclass. Let $(i,j) \in \ce$, and consider be the elementary character $\zij$ of $U_{m}(q)$. By \cite[Proposition~5.1]{A3} (see also \cite[Theorem~2.2]{ADS}), $\zij$ is constant on the superclasses of $U_{m}(q)$, and its value on the superclass associated with a basic pair $(\cd,\vphi)$ equals
\begin{equation} \label{eq:e8}
\zij(1+\ced) = \bca q^{-t} \zij(1) \tet(r \vphi(i,j)), & \text{if $(i,j) \in \cd$,} \\ q^{-t} \zij(1) , & \text{if $(i,j) \in R(\cd)-\cd$,} \\ 0, & \text{otherwise,} \eca
\end{equation}
where $t = |\set{(k,l) \in \cd}{i \prec k \prec l \prec j}|$. (We note that $q^{-t} \zij(1) = q^{t'}$ where $t'$ is the number of $\cd$-regular entries which are directly below the entry $(i,j)$.)

Using \cite[Proposition~3.2]{AN2}, we easily deduce the following result. For simplicity of writing, for any basic subset $D \sset \Phi$, we define $$R(D) = \set{\bet \in \Phi}{\ce(\bet) \sset R(\ce(D))},$$ and observe that, for any root $\bet \in \Phi$, we have $$\ce(\bet) \sset R(\ce(D)) \iff \ce(\bet) \cap R(\ce(D)) \neq \emptyset;$$ in fact, an entry $(k,l) \in \ce$ is $D$-regular if and only if $(-l,-k)$ is also $D$-regular. Further, for any root $\alp \in \Phi$, we set $$D(\alp) = \set{(k,l) \in \ce(D)}{i \prec k \prec l \prec j}$$ where $(i,j) \in \ce^{+}(\alp)$.

\begin{proposition} \label{prop:p3}
Let $\alp \in \Phi$, and suppose that $\alp \neq 2\eps_{i}$ for $1 \leq i \leq n$ (in the case where $U \leq \sp$). Let $r \in \fqx$, and $(D',\phi')$ be a basic pair for $U$. Then, $$\xa(\zdd) = \bca q^{-t(\alp,D')} \xa(1) \tet(r \phi'(\alp)), & \text{if $\alp \in D'$,} \\ q^{-t(\alp,D')} \xa(1), & \text{if $\alp \in R(D')-D'$,} \\ 0, & \text{otherwise,} \eca$$ where  $t(\alp,D') = |D'(\alp)|$.
\end{proposition}

\begin{proof}
It is enough to observe that $\xa = (\zij)_{U}$ for $(i,j) \in \ce(\alp)$ (by \cite[Proposition~3.2]{AN2}), and thus $$\xa(\zd) = \zij(\zd) = \zij(1+\ed).$$ The result follows by \refeq{e8} because $\ed = e_{\ce(D),\vphi}$ for a (uniquely determined) map $\vpd$ (by \refp{p1}).
\end{proof}

Next, we consider the case where $U \leq \sp$ and $\alp = 2\eps_{i}$ for some $1 \leq i \leq n$. Let $(D,\phi)$ be a basic pair for $U$. Let $(k,l) \in \ce$ be the smallest entry of $\ce(D)$ (with respect to the order $\preceq$ defined in \refeq{ord}), and $\bet \in D$ be such that $(k,l) \in \ce(\bet)$. Let $s = \phi(\bet)$, $D' = D - \{\bet\}$, and $\pdd$ be the restriction of $\phi$ to $D'$. We recall from the proof of \refp{p2} that $\xa(\zd) = c_{\bet,s} \xa(\zdd)$ where $c_{\bet,s} \in \C$ is given by \refeq{e7}. In particular, for $D = \{\bet\}$, we have $\zdd = 1$, hence $\xa(z_{\bet,s}) = c_{\bet,s} \xa(1) = q^{n-i} c_{\bet,s}$. In the general situation, we get $$\xa(\zd) = \xa(1) \prod_{\bet \in D} c_{\bet,\phi(\bet)} = q^{n-i} \prod_{\bet \in D} c_{\bet,\phi(\bet)},$$ and so we obtain the following formulae.

\begin{proposition} \label{prop:p4}
Suppose that $U \leq \sp$, let $\alp = 2\eps_{i}$ for some $1 \leq i \leq n$, and let $(D',\phi')$ be a basic pair for $U$.  Moreover, let $\eta$ be the quadratic character of $\fq$, let $G(\eta,\tet)$ be the Gauss sum of $\eta$ and $\tet$. Then, $$\xa(\zdd) = \bca q^{-t(\alp,D')} \xa(1) c_{\alp,r}^{D',\phi'} \tet(r \phi'(\alp)), & \text{if $\alp \in D'$,} \\ q^{-t(\alp,D')} \xa(1)  c_{\alp,r}^{D',\phi'}, & \text{if $\alp \in R(D') - D'$,} \\ 0, & \text{otherwise,} \eca$$ where $t(\alp,D') = |D'(\alp)|$, and $$c_{\alp,r}^{D',\phi'} = q^{\frac{1}{2}(t(\alp,D')-t_{0}(\alp,D'))} G(\eta,\tet)^{t_{0}(\alp,D')} \prod_{\bet \in D'_{0}(\alp)} \eta(r\inv \phi'(\bet))$$ for $D'_{0}(\alp) = D' \cap \set{2\eps_{k}}{i < k \leq n}$ and $t_{0}(\alp,D') = |D'_{0}(\alp)|$.
\end{proposition}

As an immediate consequence of \refps{p3}{p4}, we obtain following general formula for the constant value $\xd^{D',\phi'} = \xd(\zdd)$ of the supercharacter $\xd$ on the superclass $\kdd$ (see \cite[Theorem~2.2]{ADS} for the corresponding result in the case of the unitriangular group). As in the previous proposition, given any basic subset $D \sset \Phi$, we define $$D_{0} = D \cap \set{2\eps_{i}}{1 \leq i \leq n},\quad \text{and}\quad D_{0}(\alp) = D \cap \set{2\eps_{k}}{i < k \leq n}$$ whenever $\alp = 2\eps_{i} \in \Phi$ for $1 \leq i \leq n$.

\begin{theorem} \label{thm:t6}
Let $(D,\phi)$ and $(D',\phi')$ be basic pairs for $U$. Then, $$\xd(\zdd) = \bca q^{-t(D,D')} \xd(1) c_{D,\phi}^{D',\phi'} \prod_{\alp \in D \cap D'} \tet(\phi(\alp) \phi'(\alp)), & \text{if $D \sset  R(D')$,} \\ 0, & \text{otherwise,} \eca$$ where $t(D,D') = \sum_{\alp \in D} |D'(\alp)|$, and $$c_{D,\phi}^{D',\phi'} = q^{\frac{1}{2}(t(D,D')-t_{0}(D,D'))} G(\eta,\tet)^{t_{0}(D,D')} \prod_{\alp \in D_{0}} \prod_{\bet \in D'_{0}(\alp)} \eta(\phi(\alp)\inv \phi'(\bet))$$ for $t_{0}(D,D') = \sum_{\alp \in D_{0}} |D'_{0}(\alp)|$.
\end{theorem}

We observe that, in the case where, either $U \leq \eo$, or $U \leq \oo$, the set $D_{0}$ is empty, and thus $c_{D,\phi}^{D',\phi'} = 1$. Moreover, in any case, we have $q^{-t(D,D')} \xd(1) = q^{t'(D,D')}$ where $t(D,D')$ is the number of $D'$-regular entries which are directly below the entries in $\ce^{+}(D)$.

We conclude this section with a consequence of \reft{t6}, proving that every superclass factorizes uniquely as a product (in any order) of ``elementary'' superclasses; we recall that a similar factorization is valid for supercharacters.

\begin{theorem} \label{thm:t7}
Let $(D,\phi)$ be a basic pair for $U$. Then, $$\kd = \prod_{\alp \in D} K_{\alp,\phi(\alp)}$$ where the product can be taken in any order.
\end{theorem}

\begin{proof}
We order the roots according to the total order $\preceq$ defined as follows: given $\alp, \bet \in \Phi$, let $(i,j) \in \ce^{+}(\alp)$ and $(k,l) \in \ce^{+}(\bet)$, and define $\alp \preceq \bet$ if and only if $(i,j) \preceq (k,l)$ (with respect to the order on $\ce$ as defined in \refeq{ord}). On the one hand, let $z \in \prod_{\alp \in D} K_{\alp,\phi(\alp)}$ be arbitrary, and $(D',\phi')$ be the unique basic pair with $z \in \kdd$. Suppose that $(D',\phi') \neq (D,\phi)$, and let $\alp \in \Phi$ be the smallest root in $D \cup D'$ such that, either $\alp \notin D \cap D'$, or $\phi(\alp) \neq \phi(\alp')$. We consider the elementary character $\xa$ for any $r \in \fqx$. By \refps{p3}{p4}, we have $$\xa(z) = \bca c_{D,\phi} \tet(r \phi(\alp)), & \text{if $\alp \in D$,} \\ c_{D',\phi'} \tet(r \phi'(\alp)), & \text{if $\alp \in D'$,} \eca$$ where $c_{D,\phi}, c_{D',\phi'} \in \C$ are non-zero constants depending only on roots $\gam \in D \cup D'$ with $\gam \prec \alp$; moreover, by the choice of $\alp$, we have $c_{D,\phi} = c_{D',\phi'}$. Since $\tet(s) \neq 0$ for all $s \in \fq$, we conclude that $\xa(z) \neq 0$, and thus $\alp \in R(D) \cap R(D')$ (again by \refps{p3}{p4}). Now, suppose that $\alp \in D - D'$. Then, since $\alp \in R(D')$, we have $\xa(z) = c_{D',\phi'} = c_{D,\phi}$, and thus $\tet(r \phi(\alp)) = 1$. Since $r \in \fqx$ is arbitrary, we conclude that $$q = \sum_{r \in \fq} \tet(r\phi(\alp)) = \rho(\phi(\alp))$$ where $\rho = \sum_{r \in \fq} \tet_{r}$ is the regular character of $\fq^{\;+}$. It follows that $\phi(\alp) = 0$, a contradiction. Similarly, we obtain a contradiction assuming that $\alp \in D' - D$, and thus $\alp \in D \cap D'$. Thus, we get $\tet(r \phi(\alp)) = \tet(r \phi'(\alp))$ for all $r \in \fqx$, and the argument used above shows that $\phi(\alp) = \phi'(\alp)$. This final contradiction implies that $(D',\phi') = (D,\phi)$, and thus $$\prod_{\alp \in D} K_{\alp,\phi(\alp)} \sset \kd.$$

On the other hand, for the reverse inclusion, we consider the complex group algebra $\C U$ of $U$, and the superclass sum $$\hk = \sum_{z \in K} z \in \C U$$ associated with a superclass $K \sset U$. By \cite[Corollary~2.3]{DI}, the product $\prod_{\alp \in D} \hk_{\alp,\phi(\alp)}$ is a linear combination with nonnegative integer coefficients of the superclass sums of $U$. By the above, we easily conclude that $\prod_{\alp \in D} \hk_{\alp,\phi(\alp)}$ is an integer multiple of $\hk_{D,\phi}$, and so $$\kd \sset \prod_{\alp \in D} K_{\alp,\phi(\alp)},$$ as required.
\end{proof}


\section{The supercharacter table} \label{sec:table}

\reft{t4} allows the definition of the {\it supercharacter table} of $U$ as the square (complex) matrix $\ct$ having rows and columns indexed by all the basic pairs for $U$, and where the coefficient corresponding to the basic pairs $(D,\phi)$ and $(D',\phi')$ is the constant value $\xd(\zdd)$ of the supercharacter $\xd$ on the superclass $\kdd$. Since supercharacters are orthogonal (by \cite[Theorem~5.4]{AN2}), the rows of $\ct$ are orthogonal (for the usual inner product). In fact, in what follows, we prove various orthogonality relations (which are similar to the well-known for irreducible characters and conjugacy classes).

We start by considering the convolution product of supercharacters; we recall that, for any functions $\map{\zet, \eta}{U}{\C}$, the convolution product of $\zet$ and $\eta$ is the function $\map{\zet \star \eta}{U}{\C}$ defined by $$(\xi \star \zet)(z) = \sum_{x \in U} \xi(x) \zet(zx\inv)$$ for all $z \in U$. As an example, it is well-know that
\begin{equation} \label{eq:e9}
\chi \star \chi' = \del_{\chi,\chi'} \frac{|U|}{\chi(1)} \, \chi
\end{equation}
for all $\chi, \chi' \in \irr(U)$ (see, for example, \cite[Theorem~2.13]{I}); we observe that this corresponds to the generalized orthogonality relations for irreducible characters: $$\frac{1}{|U|}\, \sum_{x \in U} \chi(x) \chi'(zx\inv) = \del_{\chi,\chi'} \frac{\chi(z)}{\chi(1)}$$ for all $\chi, \chi' \in \irr(U)$ and all $z \in U$. For supercharacters, we obtain the following similar result.

\begin{theorem}[Generalized orthogonality relations for supercharacters] \label{thm:t10}
Let $(D,\phi)$ and $(D',\phi')$ be basic pairs for $U$. Then $$\xd \star \xdd = \del_{D,D'} \del_{\phi,\phi'} \frac{|U|\, \frob{\xd}{\xd}}{\xd(1)}\, \xd.$$ In other words, we have $$\frac{1}{|U|}\, \sum_{x \in U} \xd(x) \xdd(zx\inv) = \del_{D,D'} \del_{\phi,\phi'} \frac{\frob{\xd}{\xd} \xd(z)}{\xd(1)}$$ for all $z \in U$.
\end{theorem}

\begin{proof}
By \reft{t7}, we deduce that $$\chi(1) = \frob{\chi}{\rho_{U}} = \frac{\xd(1)}{\frob{\xd}{\xd}}\; \frob{\chi}{\xd},$$ and thus $$\xd = \frac{\frob{\xd}{\xd}}{\xd(1)} \sum_{\chi \in \irr_{D,\phi}(U)} \chi(1) \chi$$ for all basic pairs $(D,\phi)$. Now, since supercharacters are orthogonal, \refeq{e9} implies that $\xd \star \xdd = 0$ for all basic pairs $(D,\phi)$ and $(D',\phi')$ with $(D,\phi) \neq (D',\phi')$. On the other hand, we obtain $$\xd \star \xd = \frac{|U|\, \frob{\xd}{\xd}^{2}}{\xd(1)^{2}} \sum_{\chi \in \irr_{D,\phi}} \chi(1) \chi = \frac{|U|\, \frob{\xd}{\xd}}{\xd(1)}\, \xd,$$ as required.
\end{proof}

In particular, we deduce that the rows of the supercharacter table $\ct$ are orthogonal (an alternative proof can be easily obtained by evaluating the Frobenius scalar product).

\begin{theorem}[First orthogonality relation for supercharacters] \label{thm:t11}
Let $(D',\phi')$ and $(D'',\phi'')$ be basic pairs for $U$. Then, $$\sum_{D,\phi} \frac{|\kd|}{|U|}\, \xdd(\zd)\,\ovl{\xi_{D'',\phi''}(\zd)} = \del_{D',D''} \del_{\phi',\phi''} \frob{\xdd}{\xdd}$$ where the sum is over all basic pairs $(D,\phi)$ for $U$.
\end{theorem}

\begin{proof}
Since $\xi_{D'',\phi''}(x\inv) = \ovl{\xi_{D'',\phi''}(x)}$, the previous theorem gives $$\sum_{x \in U} \xdd(x)\,\ovl{\xi_{D'',\phi''}(x)} = \del_{D',D''} \del_{\phi',\phi''} |U| \frob{\xdd}{\xdd},$$ and the result now follows by \reft{t4}.
\end{proof}

Another consequence of the generalized orthogonality relations is the following result.

\begin{theorem} \label{thm:t12}
The space $\scf(U)$ of superclass functions is a commutative semisimple algebra with respect to the convolution product.
\end{theorem}

\begin{proof}
By \reft{t8}, the convolution product of supercharacters is a superclass function (in fact, it is a multiple of a supercharacter), hence $\scf(U)$ is a commutative algebra. Furthermore, $\scf(U)$ has a basis of orthogonal idempotents, namely the functions $$\zet_{D,\phi} = \frac{\xd(1)}{|U|\, \frob{\xd}{\xd}}\; \xd$$ for the basic pairs $(D,\phi)$, and thus it is semisimple.  
\end{proof}

Finally, we deduce that the columns of the supercharacter table $\ct$ are also orthogonal.

\begin{theorem}[Second orthogonality relation for supercharacters] \label{thm:t13}
Let $(D',\phi')$ and $(D'',\phi'')$ be basic pairs for $U$. Then, $$\sum_{D,\phi} \frac{1}{\frob{\xd}{\xd}}\, \xd(\zdd)\,\ovl{\xd(z_{D'',\phi''})} = \del_{D',D''} \del_{\phi',\phi''} \frac{|U|}{|\kdd|}$$ where the sum is over all basic pairs $(D,\phi)$ for $U$.
\end{theorem}

\begin{proof}
For any basic pairs $(D,\phi)$ and $(D',\phi')$, let $$h^{D,\phi}_{D',\phi'} = \sqrt{\frac{|\kd|}{|U|\,\frob{\xdd}{\xdd}}}\; \xdd(\zd),$$ and consider the square matrix $H = (h^{D,\phi}_{D',\phi'})$ with rows and columns indexed by the basic pairs. By the previous theorem, we have $U \bar{U}^{T} = I$, and thus $\bar{U}^{T} U = I$. It follows that $$\sum_{D,\phi} h^{D',\phi'}_{D,\phi} \ovl{h^{D'',\phi''}_{D,\phi}} = \del_{D',D''} \del_{\phi',\phi''},$$ as required.
\end{proof}


\end{document}